\documentstyle[amscd,amssymb,verbatim,12pt]{amsart}
\newcommand{\C}{{\Bbb C}}

\newcommand{\Id}{\mbox{\rm Id}}

\newcommand{\Ind}{\mbox{\rm Ind}}

\newcommand{\Int}{\mbox{\rm int}}

\newcommand{\Ker}{\mbox{\rm Ker}}

\newcommand{\R}{{\Bbb R}}

\newcommand{\supp}{\mbox{\rm supp}}

\newcommand{\vol}{\mbox{\rm vol}}

\theoremstyle{plain}
\newtheorem{definition}{Definition}

\newtheorem{theorem}{Theorem}
\newtheorem{proposition}{Proposition}
\newtheorem{corollary}{Corollary}

\numberwithin{equation}{section}
\renewcommand{\rm}{\normalshape}
\begin{document}
\title{The Dirac Operator and Conformal Compactification}
\author{John Lott}
\address{Department of Mathematics\\
University of Michigan\\
Ann Arbor, MI  48109-1109\\
USA}
\email{lott@@math.lsa.umich.edu}
\thanks{Research supported by NSF grant DMS-9704633}
\date{March 22, 2000}
\maketitle
\begin{abstract}
We give results about the $L^2$-kernel
and the spectrum of the Dirac operator
on a complete Riemannian manifold which is conformally equivalent
to the interior of a Riemannian manifold with nonempty boundary.
\end{abstract}

\section{Introduction} \label{sect1}

A general problem in spectral geometry is to understand the spectral
properties of Dirac-type operators on complete Riemannian manifolds.
In this paper we show how to use the
conformal covariance of the Dirac operator,
along with some simple arguments, to derive more precise results than are
known for general Dirac-type operators.

For background information about Dirac 
operators, we refer to 
\cite{Roe (1998)}.
Let $M$ be a connected smooth spin manifold of dimension $n \: > \: 1$. 
Let $S$ denote the spin bundle of $M$, equipped with its natural
Euclidean inner product.

Let $g$ be a complete Riemannian metric on $M$, with Dirac operator
$D_g$. Let $\Ker(D_g)$ denote the kernel of $D_g$ when acting on
$L^2(S, d\vol_g)$. 
Given $\sigma \in C^\infty(M)$, let $h$ be the Riemannian metric on 
$M$ given by 
\begin{equation} \label{conf}
g \: =\: e^{2\sigma} \: h.
\end{equation}

\begin{definition} \label{defn1}
$(M, g)$ has a conformal boundary component if one can find $\sigma$
as above and a manifold-with-boundary $N$ such that\\
1. $M$ is diffeomorphic to $\Int(N)$, \\
2. $\partial N \: \neq \: \emptyset$, \\
3. $h$ extends to a smooth Riemannian metric on $N$ and \\
4. $e^{- \: \sigma}$ extends to a locally Lipschitz function on $N$.
\end{definition}

If $(M, g)$ has a conformal boundary component then $(M, g)$ is conformally
hyperbolic in the sense of \cite{Zorich (1996)}. A basic example of
a manifold with a conformal boundary component is the real hyperbolic space
$H^n(\R)$, in which case we can take $N \: = \: 
\overline{B^n}$ and $e^\sigma(x) \: = \: 
\frac{2}{1 \: - \: |x|^2}$ for $x \in B^n$.
In general, we do not require that $N$ be compact.

\begin{theorem} \label{thm1}
If $(M, g)$ has a conformal boundary component then $\Ker(D_g) \: = \: 0$.
\end{theorem}

\begin{corollary} \label{cor1}
If $(M, g)$ has a conformal boundary component and $4 | \dim(M)$
then zero lies in the essential spectrum of $D_g$.
\end{corollary}

\begin{corollary} \label{cor2}
Let $Z$ be a closed connected spin manifold with $\widehat{A}(Z) \neq 0$.
Let $\Gamma$ be a countably infinite discrete group and let $\widehat{Z}$ be a 
connected normal $\Gamma$-cover of $Z$. If $g$ is a Riemannian metric on
$Z$, let $\widehat{g}$ be the pullback metric on $\widehat{Z}$.
Then $(\widehat{Z}, \widehat{g})$ does not admit a conformal boundary 
component. 
\end{corollary}
\noindent
{\bf Examples : } \\
1. Let $M$ be $H^n(\R)$. 
Then $M$ satisfies the conditions of 
Theorem \ref{thm1} and Corollary \ref{cor1}. In this case, the conclusions
of Theorem \ref{thm1} and Corollary \ref{cor1} were previously known by
explicit calculation \cite{Bunke (1991),Goette-Semmelmann (1999)}. \\ 
2. Let $M$ be the complex hyperbolic space $H^{2n}(\C) \: = \: G/K$, where
$G \: = \: SU(2n,1)$ and $K \: = \: S(U(2n) \times U(1))$. 
There is a Lie algebra representation of
$s(u(2n) \oplus u(1))$ on $R \: = \: \Lambda^*(\C^{2n})$ given by 
$(m, a) \rightarrow
\Lambda^*(m) \: + \: \frac{a}{2} \: \Id$. This integrates to a spinor
representation of a double cover $\widehat{K}$ of $K$ on $R$. 
Using the isomorphism $\pi_1(G) \: \cong \: \pi_1(K)$, let 
$\widehat{G}$ be the corresponding double cover of $G$. Then
$H^{2n}(\C) \: = \: \widehat{G}/\widehat{K}$ 
has a spinor bundle given by $S \: = \: \widehat{G} 
\times_{\widehat{K}} R$ and a Dirac operator $D_g$. 
In this case, $\dim(\Ker(D_g)) \: = \:
\infty$
\cite{Goette-Semmelmann (1999)}.
This shows the necessity of the condition in Theorem \ref{thm1}
that $h$ be nondegenerate on all of $N$;
note that the
numerator in the Poincar\'e metric
\begin{equation}
g \: = \: \frac{(1 - |z|^2) \sum_{i=1}^{2n} dz_i \otimes d\overline{z}_i
+ \sum_{i,j=1}^{2n} 
\overline{z}_i dz_i \otimes z_j d\overline{z}_j}{(1 - |z|^2)^2}
\end{equation}
degenerates on the boundary of the unit disk.\\
3. As remarked in \cite[Pf. of Thm. 6]{Bohr (2000)}, if $m$ is odd
then the manifold $M_{n,m}$ constructed in
\cite{Kodaira (1967)} is an aspherical spin manifold with positive signature, 
which is the total space of a surface bundle over a surface. 
Let $\widehat{M}_{n,m}$ be the universal cover; it is diffeomorphic to
$\R^4$. As $\widehat{A}(M_{n,m}) 
\neq 0$, Corollary \ref{cor2} implies that for any Riemannian metric
$g$ on $M_{n,m}$, $(\widehat{M}_{n,m}, \widehat{g})$ does not admit a conformal
compactification (in the sense of Definition \ref{defn1}) to a $4$-disk.  \\
4. The signature operator $d \: + \: d^*$ on $H^{2n}(\R)$
has an infinite-dimensional
$L^2$-kernel \cite{Donnelly (1981)}. This shows that the analog of
Theorem \ref{thm1} for general Dirac-type operators is false. \\

Now let $(Z, h)$ be a connected closed Riemannian
spin manifold.  Let $X$ be a 
closed subset of $Z$. We consider complete Riemannian manifolds $(M, g)$ which
are conformally equivalent to $(Z - X, h)$. For example,
$S^{n-m-1} \times H^{m+1}$
is conformally equivalent to $S^n - S^m$, 
where the metric on $S^{n-m-1} \times H^{m+1}$ is a product of
constant curvature metrics.

\begin{theorem} \label{thm2}
Let $(M, g)$ be a complete connected Riemannian spin manifold
of dimension $n \: > \: 1$.  Suppose that there is a 
$\sigma \in C^\infty(M)$ such
that \\
1. $(M, e^{- \: 2 \sigma} g)$ is isometrically spin-diffeomorphic 
to $(Z - X, h)$, \\
2. $\sigma$ is bounded below on $M$ and \\
3. $X$ has finite $(n-2)$-dimensional Hausdorff mass.

Then $\dim(\Ker(D_g)) \: < \: \infty$.
\end{theorem}

\begin{corollary} \label{cor3}
In addition to the hypotheses of Theorem \ref{thm2}, suppose that
$\int_M e^\sigma \: d\vol_h \: < \: \infty$. Then
$\Ker(D_g) \cong \Ker(D_h)$.
\end{corollary}

\begin{corollary} \label{cor4}
Under the hypotheses of Corollary \ref{cor3}, if $n$ is even then
the $L^2$-index of
$D_g$ is $\int_M \widehat{A}(M, g)$.
\end{corollary}

\begin{corollary} \label{cor5}
Let $(Z, h)$ be a connected closed even-dimensional Riemannian
spin manifold. 
Let $\rho$ be a nonnegative smooth function on $Z$ whose zero set is
a submanifold $X \subset Z$, along which the Hessian of $\rho$
is nondegenerate on the normal bundle $T_XZ/T_XX$.
Put $M = Z - X$, with the Riemannian metric $g \: = \: \rho^{-1} \: h$.
Then the $L^2$-index of $D_g$ is
\begin{equation}
\Ind_{L^2}(D_g) \: = \:
\begin{cases}
0 &\text{ if $\dim(X) \: = \: \dim(Z) \: - \: 1$,} \\
\int_M \widehat{A}(M, g) & \text{  if $\dim(X) \: < \: \dim(Z) \: - \: 1$.}
\end{cases}
\end{equation} 
\end{corollary}

In Section \ref{sect2} we prove the results stated in the introduction.
In Section \ref{sect3} we make some remarks.

The conformal covariance of the Dirac operator was applied to a 
different but related problem in \cite{Nistor (1999)}. I thank Victor Nistor
for pointing this out to me, and for comments on this paper.
I thank the Max-Planck-Institut-Bonn for its hospitality while this
research was performed.

\section{Proofs} \label{sect2}

If $h$ is as in (\ref{conf}),
let $D_h$ denote the Dirac operator associated to $h$. As explained in
\cite[Section 2]{Lott (1997)}, the notion of a spinor field on $M$ is
independent of the choice of Riemannian metric and so it makes sense
to compare $D_g$ and $D_h$.
From
\cite[Proposition 1.3]{Hitchin (1974)} and \cite[Proposition 2]{Lott (1986)},
\begin{equation} \label{confDirac}
D_g \: = \: e^{-\frac{(n+1) \sigma}{2}} \: D_h \:  
e^{\frac{(n-1) \sigma}{2}}.
\end{equation}

Let $\Ker(D_h)$ denote the kernel of $D_h$ on $L^2(S, d\vol_h)$.
Let $\Ker_{(\sigma)}(D_h)$ denote the kernel of $D_h$ on $L^2(S, e^\sigma \:
d\vol_h)$. 
Let ${\cal M}_{e^{\frac{(n-1)\sigma}{2}}}$ denote the operator of
multiplication by $e^{\frac{(n-1)\sigma}{2}}$ on sections of $S$.

\begin{proposition} \label{prop1}
${\cal M}_{e^{\frac{(n-1)\sigma}{2}}}$  
is an isometric isomorphism between the Hilbert spaces
$L^2(S, d\vol_g)$ and $L^2(S, e^\sigma \: d\vol_h)$ which restricts to
an isometric isomorphism between
$\Ker(D_g)$ and $\Ker_{(\sigma)}(D_h)$.
\end{proposition}
\begin{pf}
Given $\psi \in L^2(S, d\vol_g)$, put $\psi_\sigma \: = \:
e^{\frac{(n-1)\sigma}{2}} \: \psi$. Then
\begin{equation}
\int_M |\psi|^2 \: d\vol_g \: = \: \int_M |
e^{- \: \frac{(n-1)\sigma}{2}} \: \psi_\sigma|^2 \: 
e^{n \sigma} \: d\vol_h \: =
\int_M 
e^{\sigma} \: |\psi_\sigma|^2 \: d\vol_h.
\end{equation}
From (\ref{confDirac}),
\begin{equation}
D_g \: \psi \: = \: \left(
e^{-\frac{(n+1) \sigma}{2}} \: D_h \:  e^{\frac{(n-1) \sigma}{2}} \right)
\left( e^{- \: \frac{(n-1)\sigma}{2}} \: \psi_\sigma \right) \: = \: 
e^{-\frac{(n+1) \sigma}{2}} \: D_h \: \psi_\sigma.
\end{equation} 
The proposition follows.
\end{pf}
{\bf Proof of Theorem \ref{thm1} : }

Put $N^\prime \: = \: N \cup_{\partial N} ((-1, 0] \times \partial N)$, the
addition of a collar neighborhood to $N$. By assumption, $h$ extends to
a smooth Riemannian metric $h^\prime$ on $N^\prime$. Let $D_{h^\prime}$ be
the corresponding Dirac operator on $N^\prime$.

Given $\psi \in \Ker(D_g)$, put $\psi_\sigma \: =
e^{\frac{(n-1)\sigma}{2}} \: \psi$. From Proposition \ref{prop1},
$\psi_\sigma \in \Ker_{(\sigma)}(D_h)$. Let $\psi_\sigma^\prime$ be the
extension of $\psi_\sigma$ by zero to $N^\prime$. We claim that 
$\psi_\sigma^\prime$ is a smooth solution to $D_{h^\prime} \: 
\psi_\sigma^\prime \: = \: 0$. Once we show this, it will follow from
the unique continuation property of $D_{h^\prime}^2$ 
\cite[Theorem on p. 235 and Remark 3 on p. 248]{Aronszajn (1957)}, along
with the vanishing of $\psi_\sigma^\prime$ on $(-1, 0) \times \partial N$,
that $\psi_\sigma^\prime \: = \: 0$ and hence $\psi \: =\: 0$.

To show that $\psi_\sigma^\prime$ is a smooth solution to $D_{h^\prime} \: 
\psi_\sigma^\prime \: = \: 0$, by elliptic regularity theory it is enough
to show that it is a weak solution to $D_{h^\prime} \: 
\psi_\sigma^\prime \: = \: 0$. Let $\eta_\sigma^\prime$ be a smooth 
compactly-supported spinor field on $N^\prime$. Let $\eta_\sigma$ be the
restriction of $\eta_\sigma^\prime$ to $\Int(N) \cong M$. Then
\begin{equation}
\int_{N^\prime} \langle D_{h^\prime} \eta_\sigma^\prime, 
\psi_\sigma^\prime \rangle
\: d\vol_{h^\prime} \: = \: \int_M \langle D_h \eta_\sigma, 
\psi_\sigma \rangle
\: d\vol_h.
\end{equation}
We want to show that this vanishes for each choice of $\eta_\sigma^\prime$.

If $\supp(\eta_\sigma^\prime) 
\subset \Int(N)$ then we can integrate by parts to
obtain 
\begin{equation}
\int_M \langle D_h \eta_\sigma, 
\psi_\sigma \rangle \: d\vol_h \: = \: \int_M \langle \eta_\sigma, 
D_h \psi_\sigma \rangle \: d\vol_h \: = \: 0.
\end{equation}
Thus we may assume that
$\supp(\eta_\sigma^\prime) \cap \partial N \: \ne \: \emptyset$.
 
Let $K$ be a compact codimension-zero submanifold-with-boundary 
of $\partial N$. For $\epsilon$ a small
enough positive number, let $f : [0, \epsilon] \times K
\rightarrow N$ be an embedding given by Fermi coordinates near $K$.
That is, for any $k \in K$, the curve $t \rightarrow f(t, k)$ is a 
unit-speed geodesic with $f(0, k) \: = \: k$ and $(\partial_t f)(0,k) \perp
T_k \partial N$. In terms of these coordinates, we can write
\begin{equation}
h \: = \: dt^2 \: + \: m_t,
\end{equation}
where  $m_t$ is a Riemannian metric
on $K$ which depends smoothly on $t \in [0, \epsilon]$. 
To prove that $\psi_\sigma^\prime$ is a weak solution to 
$D_{h^\prime} \: \psi_\sigma^\prime \: = \: 0$
we may assume without loss of generality
that $\supp(\eta_\sigma) \subset (0, \epsilon] \times K$ for some such
$K$ and $\epsilon$. (We are thinking of $\eta_\sigma$ as being defined
on $M \cong \Int(N)$. So $\supp(\eta_\sigma)$ is a subset of $M$ which
is closed in the topology of $M$, or equivalently, in the relative topology
induced from $N$. In this sense, $(0, \epsilon] \times K$ is also a closed
subset of $M$.)

As $M$ is complete, $e^{- \: \sigma}$ vanishes on $\partial N$. Then
as $e^{- \: \sigma}$ is locally Lipschitz,
there is some $C > 0$ such that when restricted to
$f \left( [0, \epsilon] \times K \right)$,
\begin{equation}
e^{- \: \sigma} (t, k) \: \le \: C \: t.
\end{equation} 
For $t \in (0, \epsilon]$, put $K_t \: = \: f(\{t\} \times K)$.
As $\psi_\sigma \in L^2(S, e^\sigma \: d\vol_h)$, we have
\begin{equation}
\int_0^\epsilon C^{-1} \: t^{-1} \int_{K_t} |\psi_\sigma|^2 \: d\vol_{m_t}
 \: dt \: \le
\: \int_0^\epsilon \int_{K_t} e^{\sigma} \: |\psi_\sigma|^2 \: 
d\vol_{m_t} \: dt \: \le \: \int_M e^{\sigma} \: |\psi_\sigma|^2 \: 
d\vol_h \:
< \: \infty.
\end{equation}
Thus there is a sequence $t_i \in (0, \epsilon]$ such that
$\lim_{i \rightarrow \infty} t_i = 0$ and 
\begin{equation}
\lim_{i \rightarrow \infty} \int_{K_{t_i}} |\psi_\sigma|^2 \: 
d\vol_{m_{t_i}} \: = \: 0.
\end{equation}

As $\psi_\sigma \in L^2(S, e^\sigma \: d\vol_h)$, it follows that the
restriction of $\psi_\sigma$ to $f((0, \epsilon] \times K)$ is
square-integrable with respect to $d\vol_h$. Then as
$D_h \: \eta_\sigma \in L^2(S, d\vol_h)$, it follows that
$\langle D_h \: \eta_\sigma, \psi_\sigma \rangle \in L^1(M, d\vol_h)$.
Given $t \in (0, \epsilon]$, integration by parts gives
\begin{equation}
\int_{f([t, \epsilon] \times K)} 
\: \langle D_h \:  \eta_\sigma, \psi_\sigma \rangle
 \: d\vol_h \: = \: \int_{K_t} \langle 
c(\partial_t) \: \eta_\sigma, \psi_\sigma \rangle
 \: d\vol_{m_t},
\end{equation}
where $c(\partial_t)$ denotes Clifford multiplication by the unit vector
$\partial_t$.
Then 
\begin{align}
\int_M 
\: \langle D_h \:  \eta_\sigma, \psi_\sigma \rangle
 \: d\vol_h \: & = \: \lim_{i \rightarrow \infty}
\int_{f([t_i, \epsilon] \times K)} 
\: \langle D_h \:  \eta_\sigma, \psi_\sigma \rangle
 \: d\vol_h \notag \\
& = \: \lim_{i \rightarrow \infty}
\int_{K_{t_i}} \langle 
c(\partial_t) \: \eta_\sigma, \psi_\sigma \rangle
 \: d\vol_{m_{t_i}}.
\end{align}
Thus
\begin{align}
\left| \int_M 
\: \langle D_h \:  \eta_\sigma, \psi_\sigma \rangle
 \: d\vol_h \right| \: & \le \:
\lim_{i \rightarrow \infty}
\left| \int_{K_{t_i}} \langle 
c(\partial_t) \: \eta_\sigma, \psi_\sigma \rangle
 \: d\vol_{m_{t_i}} \right| \notag \\
& \le \:
\lim_{i \rightarrow \infty}
\left( \int_{K_{t_i}} |\eta_\sigma|^2
 \: d\vol_{m_{t_i}} \right)^{\frac{1}{2}} 
\left( \int_{K_{t_i}} |\psi_\sigma|^2
 \: d\vol_{m_{t_i}} \right)^{\frac{1}{2}}  \notag \\
& \le \:
\parallel \eta_\sigma \parallel_\infty \: \sup_{t \in [0, \epsilon]}
\vol^{\frac{1}{2}}(K_t) \: \lim_{i \rightarrow \infty}
\left( \int_{K_{t_i}} |\psi_\sigma|^2
 \: d\vol_{m_{t_i}} \right)^{\frac{1}{2}} \: = \: 0.
\end{align}
This proves the theorem.
\qed \\ \\
{\bf Proof of Corollary \ref{cor1} : }

Suppose that zero does not lie in the essential 
spectrum of $D_g$. Then $D_g$ is a Fredholm operator. Let $X$ be a 
connected closed $n$-dimensional 
spin manifold with $\widehat{A}(X) \neq 0$. Let $M^\prime$
be the connected sum of $M$ with $X$. Let $K \subset M$ and $K^\prime
\subset M^\prime$ be sufficiently large compact sets and let 
$g^\prime$ be a Riemannian
metric on $M^\prime$ for which $M^\prime - K^\prime$ is isometric to
$M - K$. Then $D_{g^\prime}$ is also Fredholm.
As Theorem \ref{thm1} applies to $M^\prime$, we deduce that
both $D_g$ and $D_{g^\prime}$ have vanishing $L^2$-index. However,
by the relative index theorem \cite{Gromov-Lawson (1983),Roe (1991)}, 
the difference
of the $L^2$-indices is $\widehat{A}(X)$, which is a contradiction.
\qed \\ \\
{\bf Proof of Corollary \ref{cor2} : }

Suppose that $(\widehat{Z}, \widehat{g})$ admits a conformal boundary 
component. By Theorem \ref{thm1}, $\Ker(D_{\widehat{g}}) = 0$. However, by
Atiyah's $L^2$-index theorem \cite{Atiyah (1976)}, 
$\dim(\Ker(D_{\widehat{g}})) = \infty$, which is a contradiction. 
\qed \\ \\
{\bf Proof of Theorem \ref{thm2} : }

If $\psi \in \Ker(D_g)$, put $\psi_\sigma \: = \: 
e^{\frac{(n-1)\sigma}{2}} \: \psi$. Then $\psi_\sigma$ lies in the
kernel of the Dirac operator on $Z - X$. As $\psi_\sigma \in
L^2(S, e^{\sigma} \: d\vol_h)$ and $\sigma$ is bounded below, it follows
that $\psi_\sigma \in L^2(S, d\vol_h)$. From 
\cite[p. 267, 2.3.4]{Polking (1984)}, $\psi_\sigma$ extends to an element
of $\Ker(D_h)$ on $Z$. The unique continuation property of
$D_h^2$ implies that we have constructed a
well-defined map from $\Ker(D_g)$ to $\Ker(D_h)$, which is clearly
injective.  As $Z$ is closed, $\dim(\Ker(D_h)) \: < \: \infty$.
\qed \\ \\
{\bf Proof of Corollary \ref{cor3} : }

In the proof of Theorem \ref{thm2}, we constructed maps
$\Ker(D_g) \rightarrow \Ker_{(\sigma)}(D_h) \rightarrow \Ker(D_h)$.
By Proposition \ref{prop1}, the first map is an isomorphism.
The second map is injective.  As
$\int_M e^\sigma \: d\vol_h \: < \: \infty$, the second map is surjective.
The corollary follows.
\qed \\ \\
{\bf Proof of Corollary \ref{cor4} : }

From Corollary \ref{cor3}, the $L^2$-index of $D_g$ equals the index of
$D_h$ on $Z$. From the Atiyah-Singer index theorem, this equals
$\int_Z \widehat{A}(Z, h)$. From the conformal invariance of the
$\widehat{A}$-form, this in turn equals $\int_M \widehat{A}(M, g)$.
\qed \\ \\
{\bf Proof of Corollary \ref{cor5} : }

Put $e^{2 \sigma} \: = \:
\rho^{-1}$.

If $\dim(X) \: = \: \dim(Z) \: - \: 1$, let $N$ be the metric completion
of $Z - X$ with respect to $h$. Then we can apply Theorem \ref{thm1} to 
conclude that
$\Ker(D_g) \: = \: 0$.

If $\dim(X) \: < \: \dim(Z) \: - \: 1$ then the claim follows from Corollary
\ref{cor4}. We remark that in this case, it follows from
\cite{Mazzeo (1991)} that $D_g$ has closed range, as in the model case
of $S^n - S^m$.  (When $\dim(X) \: = \: \dim(Z) \: - \: 2$, we are using
here the fact that the spin structure on $Z - X$ is inherited from the
spin structure on $Z$.) Hence from Corollary \ref{cor3}, $D_g$ is 
actually Fredholm.
\qed

\section{Remarks} \label{sect3}

Any connected
complete noncompact Riemannian manifold $(M, g)$ can be conformally
compactified to a compact metric space. For example, fix a basepoint 
$m_0 \in M$. There are  
$\phi \in C^\infty(M)$ and $c \: > \: 0$ such that
$\phi(m) \: \le \: d_g(m_0, m) \: \le \: \phi(m) \: + \: c$. Taking
a function $f : \R \rightarrow \R$ which grows sufficiently fast, 
the metric completion of $\left( M, e^{- \: f \circ \phi} \: g \right)$ will
be the Freudenthal compactification of $M$, in which one adds a point
for each end of $M$.

In general, even if $M$ has finite topological type, 
$(M, g)$ will not have a conformal compactification which is a
manifold (with or without boundary).  
The results of this paper suggest that the Dirac operator $D_g$ 
on $M$ can be studied in terms of the possibly-singular spaces which
arise as conformal compactifications of $M$. For example, let $(X, g_X)$ be a 
closed Riemannian spin manifold.
Then for any $c \: > \:0$, an end of
$M$ which is isometric to $\left( [0, \infty) \times X,
dr^2 \: + \: g_X \right)$ can be conformally compactified to the conical
space
$\left( [0,1] \times X, ds^2 \: + c^2 \: s^2 \: g_X \right)$.
In this case, the $L^2$-index theorem for manifolds with cylindrical ends 
\cite{Atiyah-Patodi-Singer (1975),Muller (1988)} 
says that there will be a contribution 
to the $L^2$-index formula for $D_g$ of $- \: \frac{1}{2}$ times the
eta-invariant of the 
link of the vertex point of the cone, i.e. of $X$. This suggests
a relationship between Dirac operators on certain complete manifolds and 
Dirac operators on singular
spaces (see \cite{Lesch (1997)} and references therein for the latter) 
although, of course, the relevant $L^2$-spaces are different.

\end{document}